# In Ratio Section Method and Algorithms for Minimizing Unimodal Functions

Vladimir Kodnyanko

Siberian Federal University, Polytechnic Institute, Krasnoyarsk, Russian Federation
E-mail: vkodnyanko@sfu-kras.ru
ORCID iD: https://orcid.org/0000-0002-2369-045X

**Abstract**. This paper proposes a new method for section an interval in a given ratio intended for minimizing unimodal functions. The ratio section search is capable of quickly recognizing monotone functions and functions with a flat bottom, which contributes to increasing its performance, as measured by the number of minimized function evaluations. The method is implemented as passive and active algorithms. A comparison of the performance of the developed method with that of the classical methods of bisection search and the golden section search was performed on the basis of the data used to minimize twenty unimodal functions of various types. For all types of functions, the passive algorithm is 2.26 times faster than the bisection search and 1.72 times faster than the golden section method. Thus, the proposed method turned out to be the fastest of the known methods of cutting off segments intended for minimizing unimodal functions. The active algorithm is faster: for all types of functions, these indicators are 3.31 and 2.52, respectively. The fastest combined Brent's method was also modernized. After the golden section procedure is replaced with a procedure for dividing a segment in a given ratio, a numerical experiment is conducted. The modernized method is 1.69 times faster than its prototype. Moreover, the performance of the active algorithm for dividing a segment at a given ratio exceeds that of the Brent's method by 1.48 times for all types of functions. The modernized Brent's method is approximately 4 times faster than the bisection search and 3 times faster than the golden section method.

## 1. Introduction

Solving various problems related to minimizing functions of one variable is an important area in the development and application of numerical mathematical methods. Using numerical methods allows many problems of not only one-dimensional but also multidimensional optimization to be solved. Success in solving them largely depends on the speed of one-dimensional optimization algorithms, so a relevant topic of research in this area is the development of methods that require the least number of calls to minimize functions.

In practical evaluations when solving optimization problems to minimize the objective minimizing functions (OMFs) of one variable, three numerical methods are most often used: uncertainty interval bisection search [1, 2, 3, 4], golden section search (GSS) [5, 6, 7, 8] and the combined Brent's algorithm [9, 10, 11, 12], which have proven their effectiveness and efficiency. The first two methods are segment elimination methods. The third method combines the reliability of the GSS and the high speed of the parabolic method. Compared with Brent's method, which is currently the fastest, the segment elimination methods have lower performance but are reliable and easy to write program code [13, 14, 15].

Unimodal functions are used in one-dimensional minimization and are often perceived differently by researchers. We adhere to the definition in [16]. According to this definition, the initial uncertainty interval $[a, b]$ can consist of three segments $[a, \alpha]$, $[\alpha, \beta]$, $[\beta, b]$ ($a \leq \alpha \leq \beta \leq b$). First, for $a < \alpha$, the function decreases; second, for $\alpha < \beta$, it is constant, i.e., it can have a flat bottom; and third, for $\alpha < b$, it increases.

For minimization methods, their reliability and speed are important and are determined by the number of function calculations. The fewer of them are, the higher the quality of the method and its efficiency. Segment separation methods have reduced speed [16, 17, 18]. However, as shown in [19], the bisection method is computationally redundant and can be improved.

## 2. Classification of minimized objective functions

All objective functions can be divided into strictly unimodal functions and others. Strictly unimodal functions always have increasing and decreasing sections and can have sections with a flat bottom ($\alpha < \beta$) or not have them ($\alpha = \beta$). Other functions can be strictly monotone—increasing or decreasing—or monotone with sections of a flat bottom on one of the segments of the definition domain. A special case of the latter is a constant function. In accordance with this classification, it is possible to construct minimization algorithms that can work faster than algorithms that do not make exceptions for nonstrict unimodal functions.

Let us consider several approaches to minimizing such functions.

1. Constant function (Table 2, cell 1/1). A rare but possible class of unimodal functions for which ($a$ = α) & (β = $b$). This function has a flat bottom throughout its domain. Obviously, it can be determined in just three OMF evaluations. For such functions, the mentioned methods work more slowly because when searching for a solution, they do not distinguish between functions by their type.

2. Monotone functions (Table 2, cells 1/2, 1/3). Theoretically, their minima can be found in just three OMF evaluations. The scheme is as follows. First, the function is evaluated at two points $a + \varepsilon_0$ and $b - \varepsilon_0$, where $\varepsilon_0$ is the tolerance (accuracy, admissible error). On this basis, it is presumably possible to establish the nature of the monotonicity of the function. If $f(a + \varepsilon_0) < f(b - \varepsilon_0)$, then the function presumably monotonically increases. Now, we need to evaluate $f(a)$. If $f(a) < f(a + \varepsilon_0)$, then the function is strictly monotonically increasing. If $f(a) = f(a + \varepsilon_0)$, then the OMF is a function with a flat bottom. In both cases, the solution to the problem is $x = a$. Otherwise, such a function is strictly unimodal. The function is checked for a monotonic decrease in a similar way. The OMF should not be checked for monotonicity after it is evaluated at the first two points because the algorithm checks almost all functions for monotonicity, which leads to a loss in the algorithm's performance. It is better to carry out this procedure after evaluating the first four or even five points. If they line up in a monotonic sequence, which is rarely possible for strictly unimodal functions, then only after this should one or two extreme points of the segment be checked. In the evaluations of monotonicity, a variant with four points was used. The algorithm determines the monotonic function and its minima in six OMF evaluations.

Table 1. Graphs of strictly unimodal functions without a flat bottom

| | 1 | 2 | 3 |
|---|---|---|---|
| 1 | | | |
| 2 | | | |

Table 2. Graphs of quickly recognizable unimodal functions

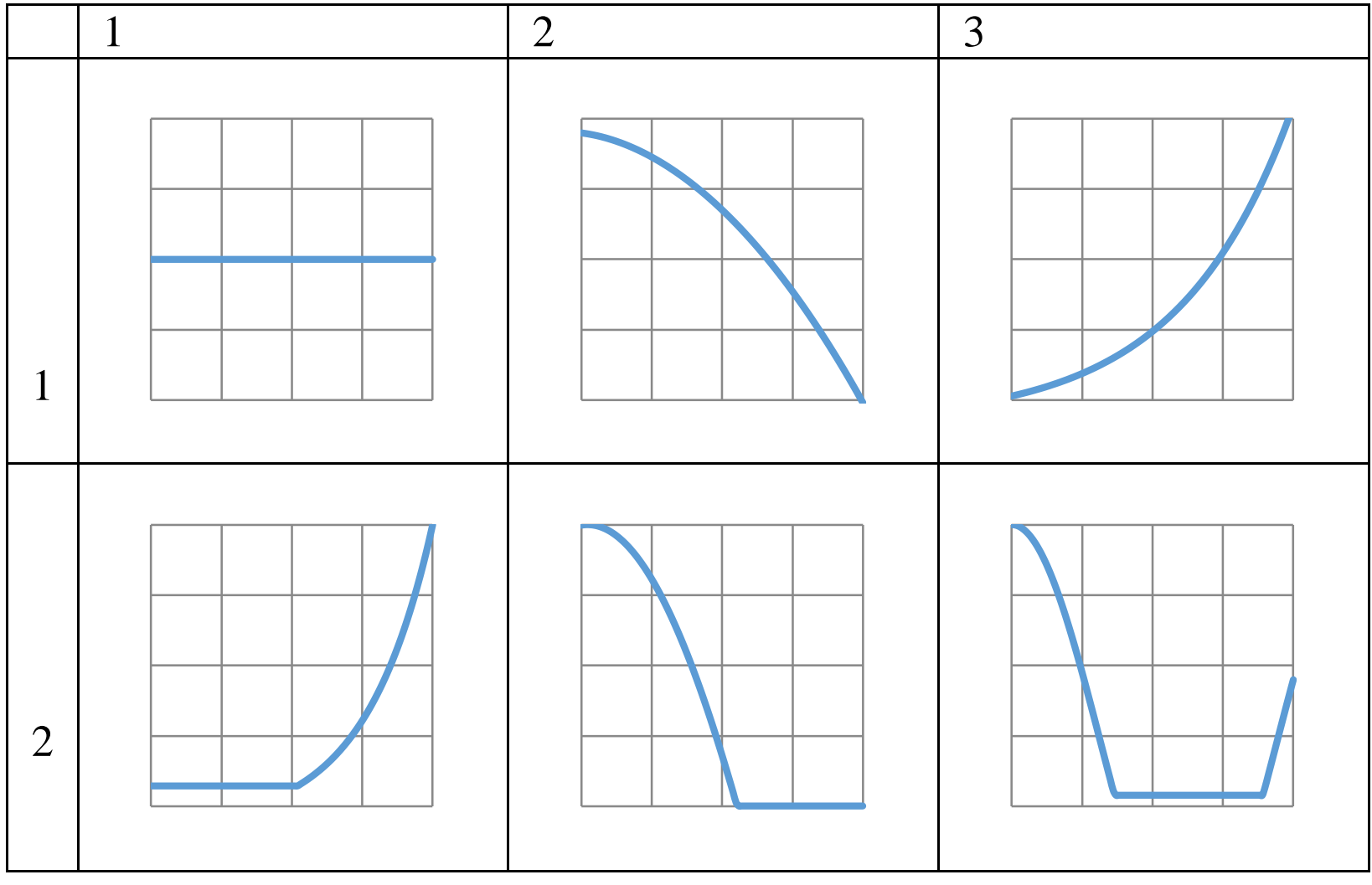

3. Functions with a flat bottom (Table 2, cells 2/1, 2/2, 2/3). An example of such a constant function has already been considered, but any other function can have a flat bottom. If any algorithm finds three points with the same ordinates, then, owing to the unimodality of the OMF, any of them is a solution to the problem.

## 3. Block diagram for recognizing and minimizing strictly monotone functions

Fig. 1 shows a block diagram of the Boolean algorithm function *Mnt* for recognizing strictly monotone functions, examples of which are shown in the graphs of Table 2, cells 1/2, 1/3.

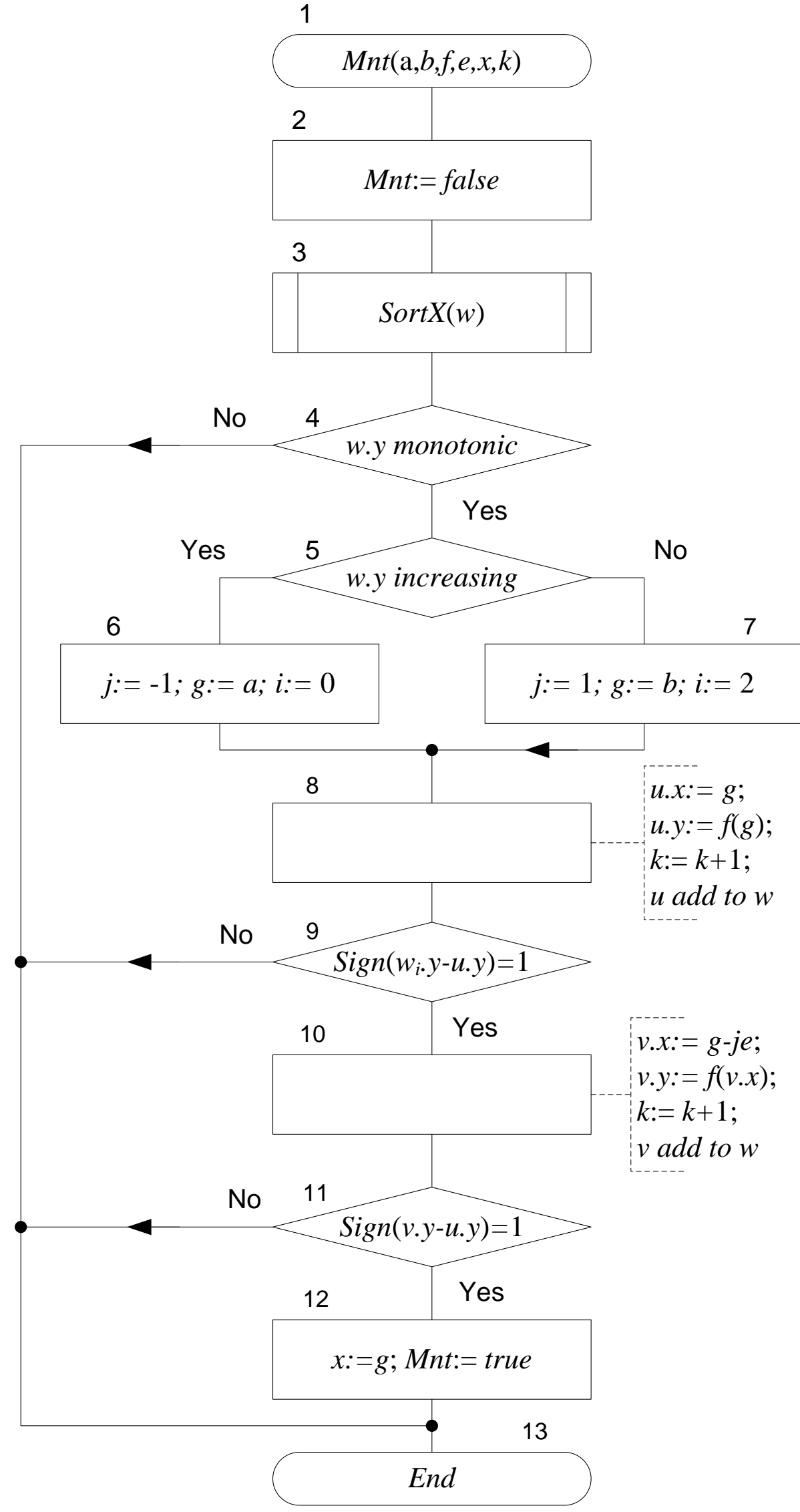

Fig. 1. Block diagram of the algorithm function *Mnt*

As mentioned, for any strictly monotonic function, the solution is achieved in six evaluations of the function (four evaluations to propose a hypothesis about the monotonicity of the function before evaluating this function and one or two more evaluations in the function body to refute or confirm the hypothesis). Two-component variables for points of the form $p = p\ (x, y)$ are used to designate points.

When compiling the algorithm for the function *Mnt* ($a$, $b$, $f$, $k$, ε, $w$), the following were used:

- boundaries $a$, $b$ of the original uncertainty interval,
- name $f$ of the minimized function,
- number of function evaluations $k$,
- tolerance ε,
- array $w$ of four two-component points that have already been evaluated in the process of minimizing the function via the algorithm function *Mnt*.

At the beginning of the algorithm in block 2, the value *Mnt* = *false* is assigned under the assumption that the proposed hypothesis about the monotonicity of the digital function in the process of the operation algorithm may be false.

Next, in block 3, the array *w is* sorted in ascending order of the abscissas of its elements. In block 4, the algorithm determines whether the sequence composed of the ordinates of the array *w* is monotonic. If this is confirmed, then depending on whether it is increasing or decreasing, which is determined in block 5, the variables *j*, *g*, and *i* in one of blocks 6 or 7 are assigned certain values depending on the nature of the monotonicity of the function. If the hypothesis is already refuted at the initial stage of the analysis of this nature, then the work of the algorithm function ends.

Then, in block 8, a new point *u* of the OMF is evaluated, which is one of the ends of the original uncertainty segment. In block 9, a check is performed to determine whether it satisfies the proposed hypothesis. If the point does not correspond to it, the process ends; otherwise, in block 10, the point *v* is computed, which is located from the required end of the segment by the tolerance *e*.

If point *v* also corresponds to the hypothesis in block 11, then in block 12, the algorithm function receives the value *Mnt* = *true*, which indicates that the minimizing function is monotonic. It would be possible to check whether three points with the same ordinates appeared in the array *w* when two new points, *u* and *v,* were added to it. If such points are present in the array, this means that in the process of recognizing the monotonicity of the function, a flat bottom is simultaneously detected. Therefore, the abscissa of any of these points, owing to the unimodality of the function, is a solution to the problem, which also indicates the completion of the minimization process. However, it is better to perform such a check immediately after the false completion of the *Mnt* algorithm.

## 4. In ratio section algorithm *RatioP*

As mentioned earlier, there is a possibility of improving the performance of the algorithm when the OMF is minimized because only one evaluation is required at each iteration, including the first iteration.

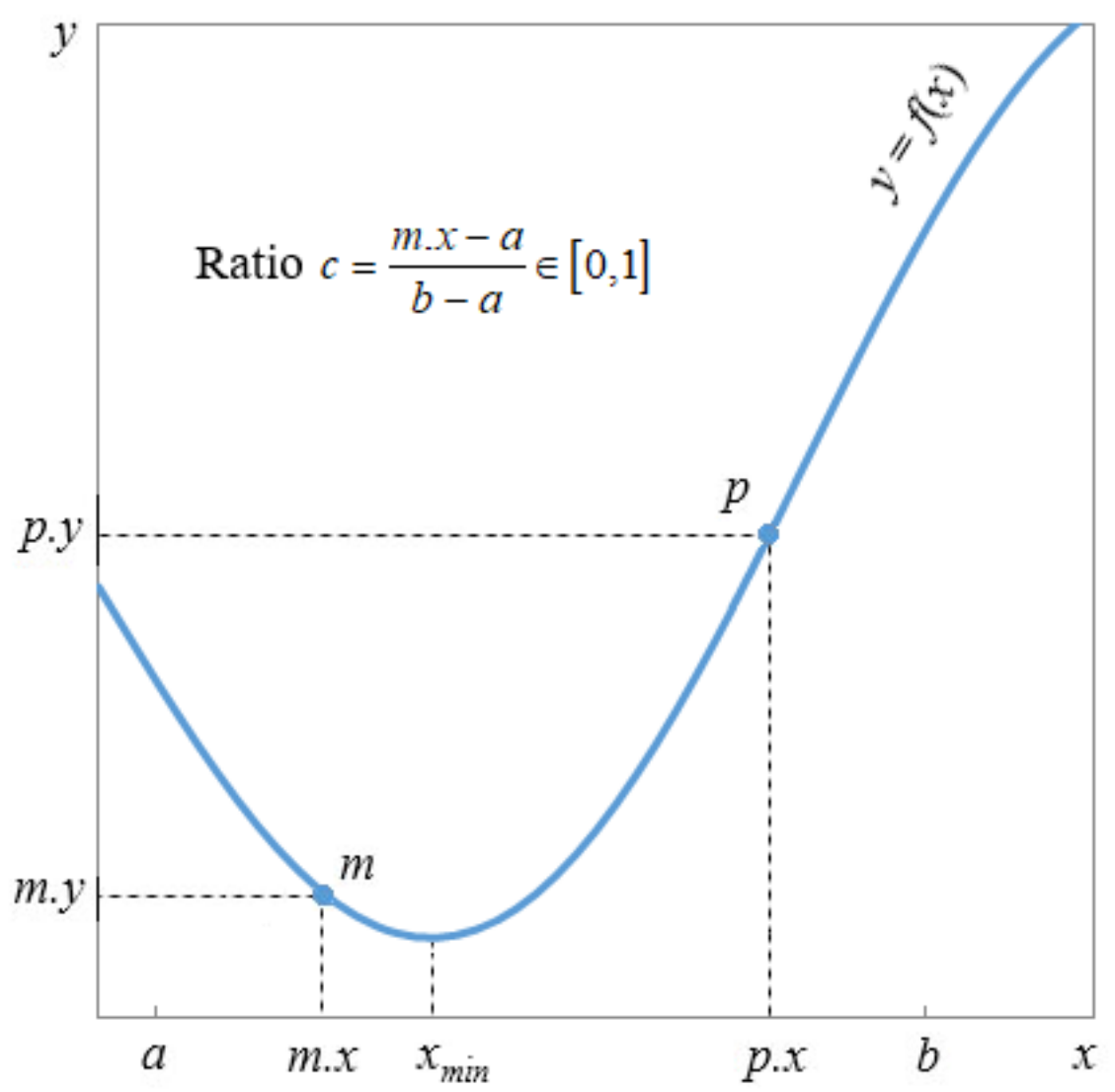

Fig. 2. Graph of minimized function *f*(*x*)

Below is the *RatioP* (*a*, *b*, *c*, *f*, *k*, *x*, *w*) algorithm with parameters similar to those of the *Mnt* function, the algorithm of which is shown in Fig. 1. Graph of minimized function *f*(*x*) is shown in Fig. 2. The constant $c = \frac{m.x - a}{b - a} \in [0,1]$ is a ratio by which a new point *p* can be set on one of the segments (*a*, *m.x*) or (*m.x*, *b*) of the uncertainty interval [*a*, *b*] closer to or further from one of its edges, as is the case when dividing a segment in a given ratio *c*, where *m.x* is a point on this segment. At *c* = 0.5, such a point is set in the middle of the segment, which corresponds to the bisection mode. At *c* < 0.5, the point will be closer to *m.x*, and at *c* > 0.5, it will be farther and closer to one of the ends of the segment (*a* or *b*).

A diagram of the algorithm is shown in Fig. 3. After the *RatioP* algorithm starts working, in block 2, a point is evaluated in the middle of the segment [*a*, *b*], which becomes the first point with the smallest ordinate. Next, the initial number of function evaluations *k* and the number of points with the smallest ordinate *n* are determined, the array *w* of evaluated points is cleared, and the first point *m* is added to this array. Then, in blocks 3 – 18, which form a cycle, the process of finding function minima with a given tolerance is launched.

At each step of the cycle, in block 3, from two segments [*a*, *m.x*] and [*m.x*, *b*] of segment [*a*, *b*], a segment of greater length is selected. In blocks 4 or 5, the abscissa *p.x* of the new point *p* is determined by taking into account the

value of the ratio $c$

$$p.x = \begin{cases} ca + (1-c)m.x, & m.x - a > b - m.x, \\ cb + (1-c)m.x, & otherwise. \end{cases}$$

In block 6, the ordinate of this point is evaluated

$$p.y = f(p.x)$$

and the point is added to the array $w$.

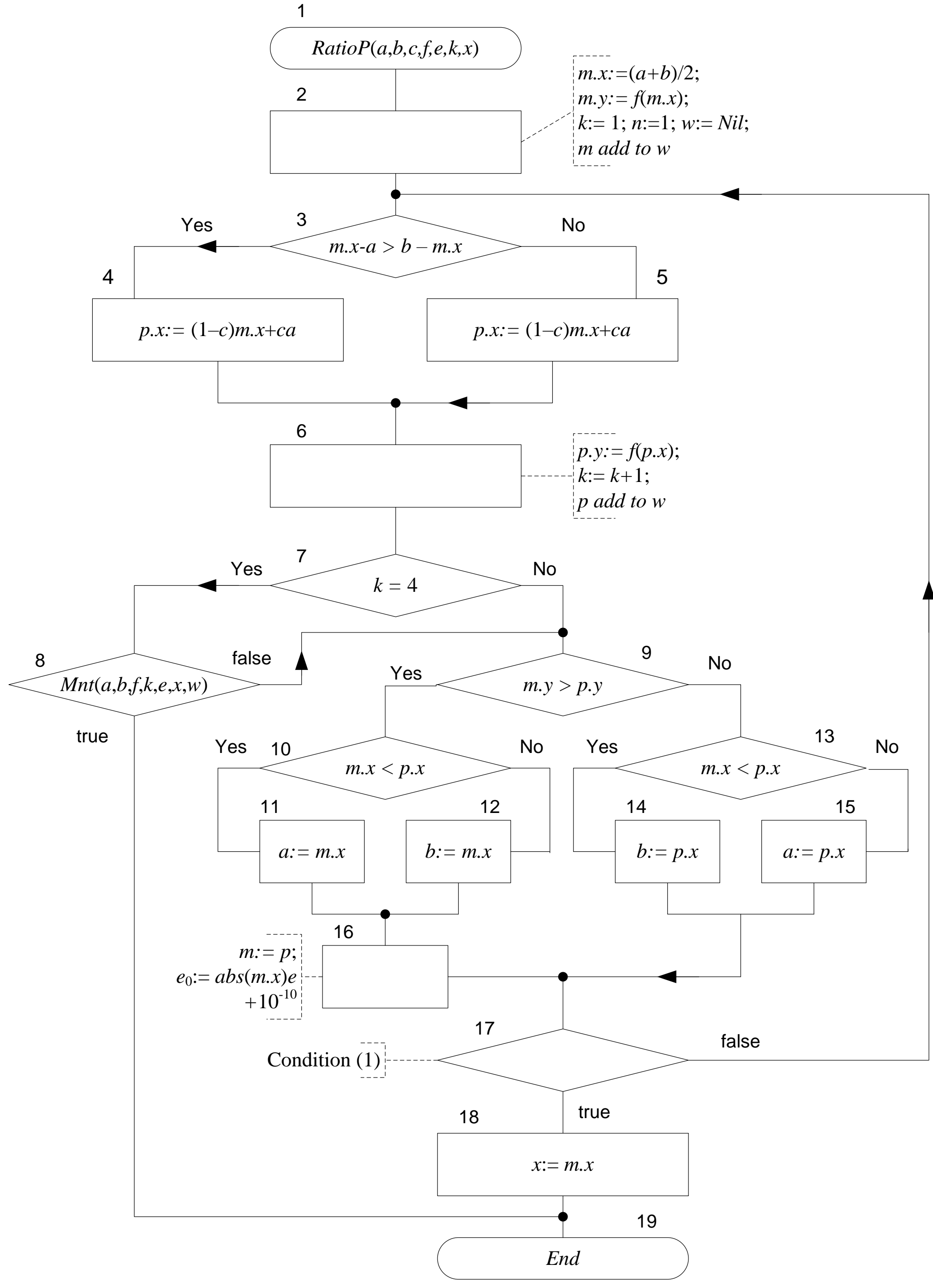

Fig. 3. Block diagram of the *RatioP* algorithm

In block 7, at the third step of the cycle in array $w$, the number of points has reached the value $k = 4$; in block 8, the algorithm *Mnt*, which is shown in Fig. 1, is evaluated to check the function for monotonicity. If the monotonicity hypothesis is confirmed, then having determined the abscissa of the minimum of the function via this algorithm, the *RatioP* algorithm finishes its work.

If at this step, the monotonicity hypothesis is not confirmed or $k \neq 4$, then in block 9, the ordinates of the new point

$p$ and the point of the current minima $m$ are compared. If condition 9 is true, this means that a new minima of the function has been found, and in blocks 11 and 12, the abscissa $m.x$ of the previous minima becomes the new boundary of the uncertainty segment depending on condition 10, which indicates on which side of $m$ the abscissa of the new point $p$ is located. In block 10, the current minimum is replaced by a new minimum.

The conditions for stopping the process are as follows:

$$\left| m.x - \frac{a+b}{2} \right| <= e_0 - \frac{1}{2}(b-a), \tag{1}$$

which used in Brent's algorithm [9], where $e_0 = e|m.x| + 10^{-10}$.

If this condition is not met, i.e., a new minima point not found at the current iteration, then in block 13, depending on the fulfillment of its condition, it will be determined which of the boundaries of the uncertainty segment should be changed in blocks 14 or 15 using point $p$.

Since at each step of the cycle, the segments are separated, the length of the original uncertainty segment decreases; then, in block 17, the condition checks whether the algorithm has achieved its goal. The process stops when one of the conditions is met, when the length of the longest segment has become less than the required tolerance, or when the $N3$ function is used for array w, at least three different points with the same ordinates are found. This indicates that at the step of the cycle, a flat bottom of the OMF graph is determined.

If none of these conditions are met, then a transition to block 3 is made to perform new steps of the cycle until the result is achieved.

## 5. Comparative results of the numerical experiment using the RatioP algorithm

To study the effectiveness of the *RatioP* algorithm, as well as to assess the practical prospects of the algorithm, a program code was developed in the Delphi language [19]. Computations were performed for various types of unimodal functions: smooth, flat, slowly changing, piecewise, monotonic, and flat bottom functions. Evaluations were performed for $\varepsilon = 10^{-5}$ via variables and constants of an arithmetic type that supports 19–20 significant digits. Golden section and bisection searches were used as control methods.

Table 3. Functions used in the studies and their uncertainty intervals

| $N$ | $f(x)$ | $[a, b]$ | $N$ | $f(x)$ | $[a, b]$ |
|---|---|---|---|---|---|
| 1 | $1$ | [0.5, 1.5] | 11 | $0.3+\cos(x^2+2x-3)$ | [–0.9, 0.9] |
| 2 | $20+\dfrac{16}{x}$ | [2.6, 6.8] | 12 | $0.2+(x-1.5)^2$ | [0.3, 3.2] |
| 3 | $1.5+Exp(x)$ | [1.2, 3.7] | 13 | $100+\left[1-\exp(x)\sin(x)\right]^2$ | [0.1, 1] |
| 4 | $1.5+\mathrm{Max}\left[4\cos(x),1\right]$ | [0.1, 4.9] | 14 | $1.2-\cos(x^2)$ | [–1.2, 1.5] |
| 5 | $1.2+\mathrm{Max}\left[5exp(x)-1,1\right]$ | [–1.6, 1.1] | 15 | $1.2+5\exp(-x^2)+x$ | [0.3, 11.5] |
| 6 | $1.5+\mathrm{Max}\left[\cos(4-x^2),0.5\right]$ | [3.2, 3.5] | 16 | $1.2+\exp(-x)+3.5\sin(x)$ | [–1.6, 0.8] |
| 7 | $1.5+\mathrm{Max}\left[\exp(-x),\cos(x),x^4,x^2\right]$ | [–0.6, 1.1] | 17 | $2.3+3\exp(x)-x^2+5x$ | [–15, 7] |
| 8 | $0.2+\mathrm{Max}\left[13(x-2)^2,20(x-1)\right]$ | [–1.2, 2.7] | 18 | $1.2+3\mathrm{ch}(x-2)-2\mathrm{sh}(x-3)$ | [–2.1, 2.5 |
| 9 | $1.2+\lvert x-1\rvert$ | [0.5, 6.5] | 19 | $2.3+\left[\exp(3-x)+4(x-2)\right]^2$ | [–0.5, 2.5] |
| 10 | $12+1000\lvert(x-2)\rvert^{8.4}$ | [1, 4.3] | 20 | $1.2+\lvert x-2\rvert^{3.6}$ | [0.1, 1] |

Twenty functions were used to process the evaluation results. The types of functions and their variation intervals are given in Table 3. Among them, 1 is a constant (1), 2 is strictly monotone (2, 3), 3 has a flat bottom (4 – 6), and 15 are other strictly unimodal (14 – 20).

Table 4 shows the data on the number of OMF evaluations via golden section search, bisection search and the proposed *RatioP* method for the ratio $c = 0.5$, which, as mentioned, exactly corresponds to bisection.

In Table 4 and the subsequent similar Tables, in each row whose number corresponds to the number of minimized OMFs, the cell corresponding to the fastest method is highlighted in green, the next fastest method is highlighted in light green, and the one after that is highlighted in gray.

Compared with the golden section search and *RatioP* algorithms, the classical bisection method (column title *Bisec*) requires a greater number of OMF evaluations for all the studied functions. The $\Sigma k$ row shows the total number of OMF evaluations that were required to minimize all twenty OMFs in Table 3. For this purpose, the bisection algorithm had 768 OMF evaluations, 585 golden section evaluations (column title *Gold*), and 465 *RatioP* evaluations. That is, the bisection search worked the slowest, and the *RatioP* algorithm worked the fastest. As follows from the *Relat* row, which contains the data on the number of OMF evaluations reduced to *RatioP*, GSS took an intermediate position among the methods under study, yielding a performance *RatioP* of 1.26, whereas the bisection search was inferior to *RatioP* by 1.65 times.

The first reason for the better performance is that the first 6 functions are recognized faster by the algorithm.

As shown in Table 4, for *RatioP,* only three OMF evaluations are required to recognize a constant, while bisection search requires 36 evaluations, and GSS requires 26 evaluations.

Table 4. Evaluations to all functions 1 – 20

| *N* | *Bisec* | *Gold* | *RatioP* |
|---|---|---|---|
| 1 | 36 | 26 | 3 |
| 2 | 32 | 25 | 6 |
| 3 | 36 | 28 | 6 |
| 4 | 36 | 27 | 4 |
| 5 | 36 | 29 | 8 |
| 6 | 28 | 21 | 4 |
| 7 | 36 | 28 | 26 |
| 8 | 38 | 29 | 31 |
| 9 | 40 | 30 | 31 |
| 10 | 34 | 27 | 27 |
| 11 | 44 | 33 | 37 |
| 12 | 36 | 28 | 31 |
| 13 | 36 | 27 | 29 |
| 14 | 70 | 51 | 33 |
| 15 | 40 | 30 | 34 |
| 16 | 38 | 29 | 29 |
| 17 | 42 | 32 | 35 |
| 18 | 38 | 29 | 31 |
| 19 | 38 | 30 | 32 |
| 20 | 38 | 28 | 30 |
| $\Sigma k$ | 772 | 587 | 467 |
| Relat | 1.65 | 1.26 | 1 |

Table 5. Evaluations of strictly unimodal functions 7 – 20

| *N* | *Bisec* | *Gold* | *RatioP* |
|---|---|---|---|
| 1 | – | – | – |
| 2 | – | – | – |
| 3 | – | – | – |
| 4 | – | – | – |
| 5 | – | – | – |
| 6 | – | – | – |
| 7 | 36 | 28 | 26 |
| 8 | 38 | 29 | 31 |
| 9 | 40 | 30 | 31 |
| 10 | 34 | 27 | 27 |
| 11 | 44 | 33 | 37 |
| 12 | 36 | 28 | 31 |
| 13 | 36 | 27 | 29 |
| 14 | 70 | 51 | 33 |
| 15 | 40 | 30 | 34 |
| 16 | 38 | 29 | 29 |
| 17 | 42 | 32 | 35 |
| 18 | 38 | 29 | 31 |
| 19 | 38 | 30 | 32 |
| 20 | 38 | 28 | 30 |
| $\Sigma k$ | 568 | 431 | 436 |
| Relat | 1.30 | 0.99 | 1 |

Monotone functions are recognized by the *RatioP* algorithm in 6 OMF evaluations, whereas bisection search requires 32–36 and GSS requires 25–28. As shown in Table 4, *RatioP* requires a comparatively smaller number of OMF evaluations to recognize functions 4–6 with a flat bottom.

If such functions, which, although rare, can be included in the evaluations during optimization, are excluded, different results will be obtained. The data for functions 7–20 are given in Table 5. The results clearly show that bisection search is still inferior to *RatioP* in terms of performance, but GSS is no longer inferior to the studied method and yields practically the same result.

This leads to the important conclusion that, similar to the GSS, only one evaluation of the OMF is required at each iteration, resulting in practically the same results. This contradicts the deeply rooted opinion that the GSS works significantly faster than the bisection method does. If we add that *RatioP* recognizes functions of the mentioned type faster, then it works faster than the GSS. Notably, while the golden section search does not have a resource for improving performance, the *RatioP* section search does. The expectation is that the performance of *RatioP* can be improved by optimizing the algorithm with respect to the ratio *c*.

## 6. Algorithm for smoothing the results of a numerical experiment

In the process of optimizing the performance of the *RatioP* algorithm by the ratio *c*, data such as any experimental results obtained from a limited number of experiments, characterized by data scatter, are obtained. To smooth the data, an algorithm based on the root-mean-square deviation method was developed [20].

A block diagram of the *LS* (least squares) algorithm procedure is shown in Fig. 4.

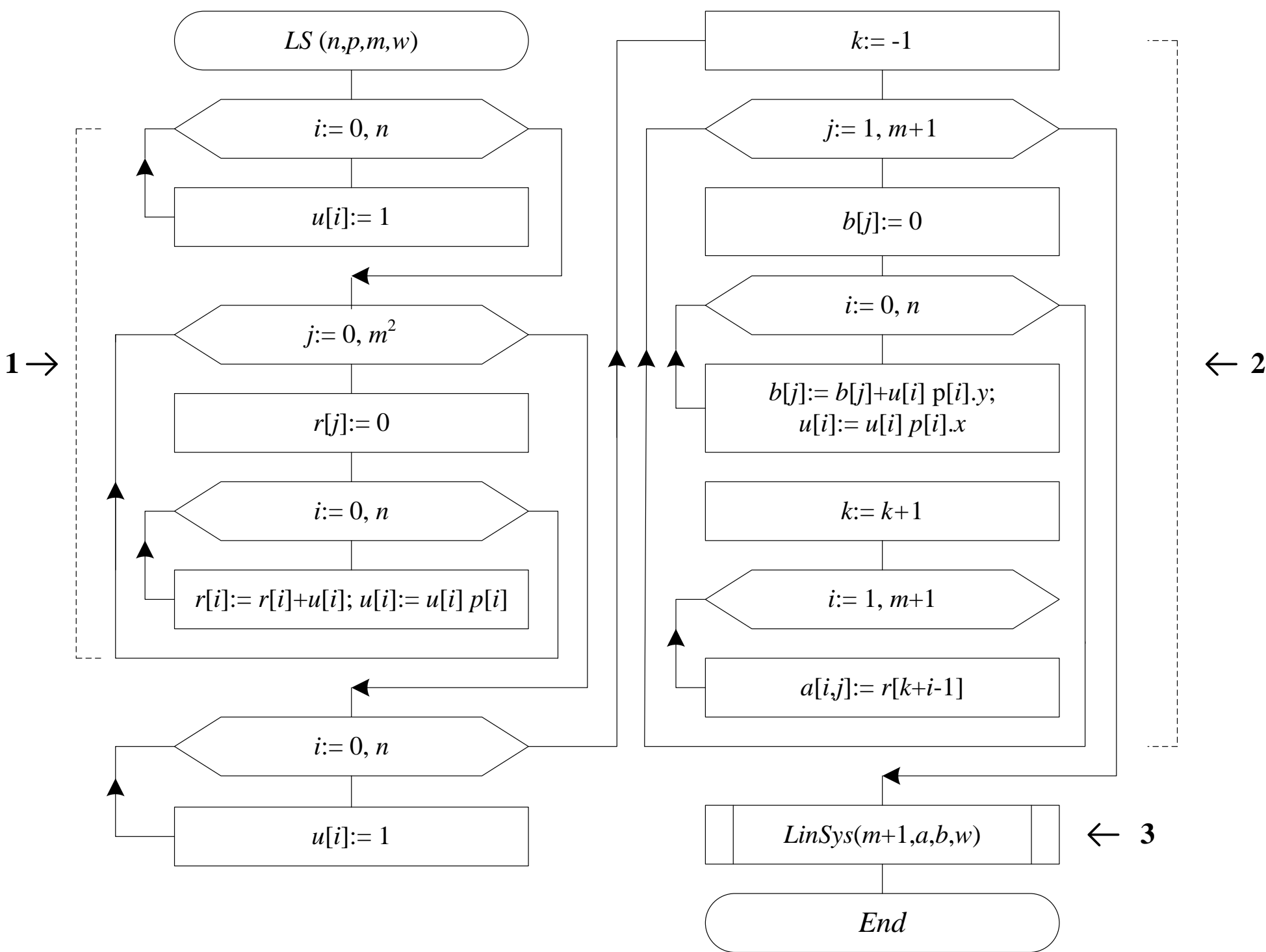

Fig. 4. Block diagram of the least squares (*LS*) algorithm

The input parameters are the quantitative characteristic of the experimental points *n*, the array of experimental points *p*[0], *p*[1], …, *p*[*n*], and the desired degree of the resulting polynomial *m*. The output data are the array of polynomial coefficients *w*[0], *w*[1], …, *w*[*m*]. In block group 1, the array *r*[0], *r*[1], …, *r*[$m^2$] is formed, which is then used by block group 2 to form the matrix *a* of the system of linear equations, the solutions to which are the coefficients of the polynomial.

The elements of the array *r*, the a matrix, and the *b* vector of the right-hand sides of the equations are evaluated in block group 2 via the formulas

$$r_j = \sum_{i=0}^{n} p_{i,x}^{j}, b_j = \sum_{i=0}^{n} p_{i,y} p_{i,x}^{j}, a_{i,j} = r_{k+i-1}, k = 0,1,...,m.$$

In block 3, the procedure for solving the system of equations accessed, the result of which is an array of coefficients *w*. To solve this system of equations, the Gauss method was used.

The algorithm has been verified by testing on arrays in the ranges $n = 3 – 100$, $m = 3 – 10$, and $m < n$.

## 7. Results of *RatioP* algorithm optimization by ratio *c*

Figure 5 shows the graph of the dependence of the average number of evaluations *K* of the *RatioP* algorithm on all the OMFs of Table 3 for different values of the ratio *c*. The value of *c* varies within the segment $c \in (0.01, 0.80)$ with a step of 0.01.

The graph shows that this dependence is extreme in nature, with clearly defined minima in the vicinity of the value $c = 0.2$. The solid line in the graph corresponds to the experimental data, and the dots represent the smoothing curve obtained via the *LS* procedure (Fig. 4). The smoothing curve with the order of the polynomial $l = 5$ has the best correspondence to the experimental data.

The graph in Fig. 5, *b* shows that the minimum average number of evaluations *K* of functions 7–20 is noticeably less than the similar value of the criterion $K$ for the ratio $c = 0.5$, corresponding to bisection. Indeed, for $K(0.2) \approx 17$, whereas when the *RatioP* algorithm operates in the bisection mode, $K(0.5) \approx 22$.

There is a simple explanation for the phenomenon discovered. If we refer to the operators of blocks 4 and 5 of the *RatioP* algorithm flowchart, which is shown in Fig. 2, we can see that the smaller $c$ *is*, the closer the new point $p$ is to the point of the current minimum. It is reasonable to assume that in comparison with the bisection parameter $c = 0.5$, the smaller $c$ is, i.e., the closer the new point is to the current minimum point, the more effectively the uncertainty segment will be compressed. However, the value of the ratio $c$ can be reduced only to a certain limit, since, for example, at $c = 0$, the new point will coincide with the current minimum point and, therefore, in this case, the algorithm will loop, depriving it of efficiency.

Tables 6 and 7 contain data similar to Tables 4 and 5 for $c = 0.2$.

As shown in Table 6, the *RatioP* algorithm, when all twenty functions are optimized with $c = 0.2$, works 1.72 times faster than the GSS search and 2.26 times faster than the bisection search. If we exclude special functions 1–6 to put the *RatioP* algorithm on an equal footing with classical methods, as shown in Fig. 2, it will still have better performance indicators, outperforming GSS in performance by 1.40 times and bisection search by 1.84 times.

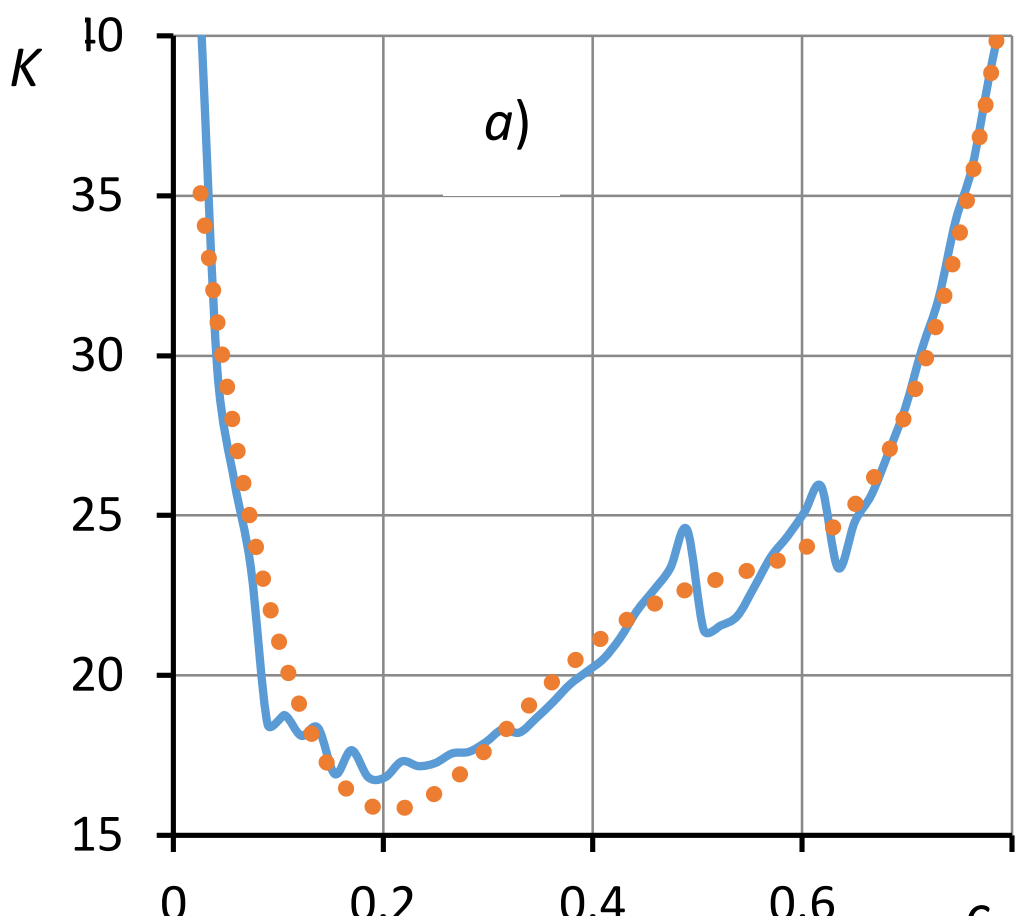

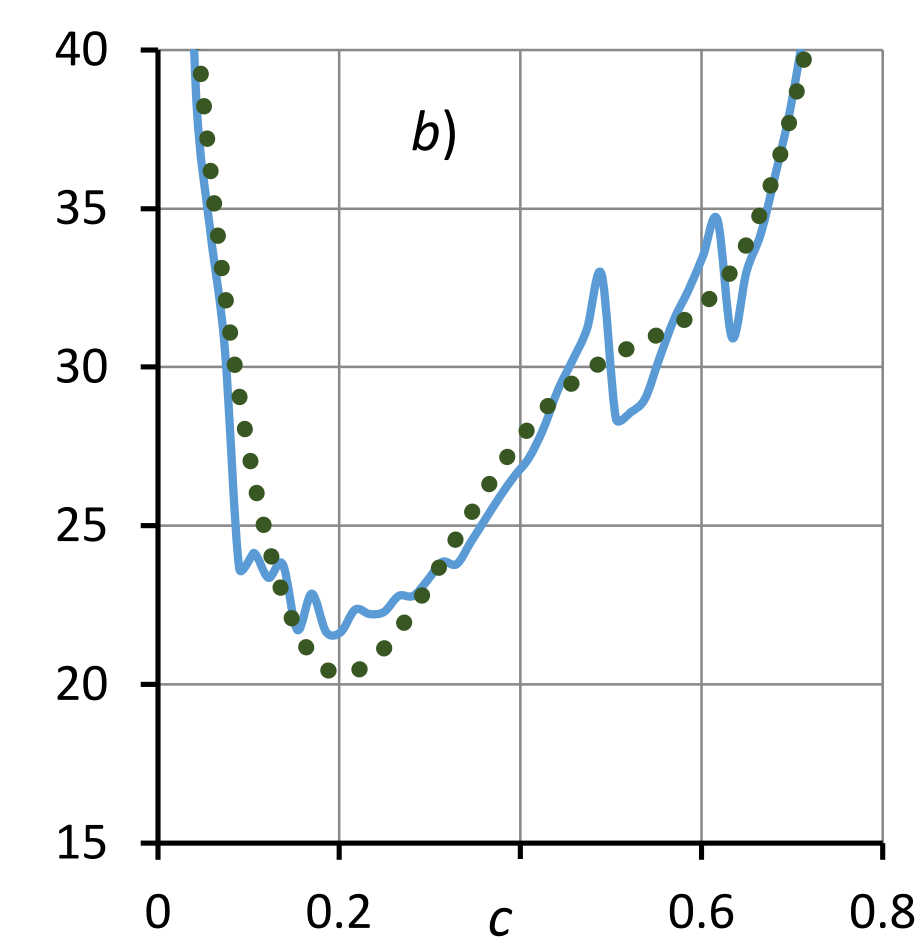

Fig. 5. Graph of the average number $K$ of evaluations according to the *RatioP* algorithm versus the ratio $c$:
*a*) – for functions 1–20, *b*) – for functions 7–20

Table 6. Evaluations of all functions 1 – 20

| *N* | *Bisec* | *Gold* | *RatioP* |
|---|---|---|---|
| 1 | 36 | 26 | 3 |
| 2 | 32 | 25 | 6 |
| 3 | 36 | 28 | 6 |
| 4 | 36 | 27 | 4 |
| 5 | 36 | 29 | 10 |
| 6 | 28 | 21 | 4 |
| 7 | 36 | 28 | 23 |
| 8 | 38 | 29 | 21 |
| 9 | 40 | 30 | 27 |
| 10 | 34 | 27 | 18 |
| 11 | 44 | 33 | 26 |
| 12 | 36 | 28 | 20 |
| 13 | 36 | 27 | 18 |
| 14 | 70 | 51 | 21 |
| 15 | 40 | 30 | 26 |
| 16 | 38 | 29 | 22 |
| 17 | 42 | 32 | 23 |
| 18 | 38 | 29 | 20 |

Table 7. Evaluations of strictly unimodal functions 7 – 20

| *N* | *Bisec* | *Gold* | *RatioP* |
|---|---|---|---|
| 1 | – | – | – |
| 2 | – | – | – |
| 3 | – | – | – |
| 4 | – | – | – |
| 5 | – | – | – |
| 6 | – | – | – |
| 7 | 36 | 28 | 23 |
| 8 | 38 | 29 | 21 |
| 9 | 40 | 30 | 27 |
| 10 | 34 | 27 | 18 |
| 11 | 44 | 33 | 26 |
| 12 | 36 | 28 | 20 |
| 13 | 36 | 27 | 18 |
| 14 | 70 | 51 | 21 |
| 15 | 40 | 30 | 26 |
| 16 | 38 | 29 | 22 |
| 17 | 42 | 32 | 23 |
| 18 | 38 | 29 | 20 |

| 19 | 38 | 30 | 22 |
|---|---|---|---|
| 20 | 38 | 28 | 21 |
| Σk | 772 | 587 | 341 |
| Relat | 2.26 | 1.72 | 1 |

| 19 | 38 | 30 | 22 |
|---|---|---|---|
| 20 | 38 | 28 | 21 |
| Σk | 568 | 431 | 308 |
| Relat | 1.84 | 1.40 | 1 |

Additionally, mass evaluations were carried out for 40 different functions and 10 variants of the uncertainty interval boundaries so that the function on the interval $[a, b]$ remained strictly unimodal and did not have a flat bottom. An evaluation of 400 such problems revealed that the ratio section search algorithm *RatioP* had better performance indicators than did the GSS by 1.47 times and the bisection search by 1.91 times. Apparently, these indicators are closer to comparative average statistical estimates.

The proposed method, even at $c = 0.5$, is not, in a certain sense, an improvement over the classical bisection method since the principle underlying the *RatioP* algorithm does not require determining the nature of the function on the increase or decrease in the OMF at a new point. It is based on the principle of the fastest reduction of the uncertainty segment, which does not require such data; i.e., it is a completely new method. To avoid confusion with the bisection method, it should be evaluated differently. In essence, this is a method of ratio division of a segment with a ratio $c$, so it can be evaluated exactly that. The results of the conducted research indicate that, on average, the *RatioP* works approximately 2 times faster than the bisection search and 1.5 times faster than the GSS.

## 8. *RatioA* Algorithm

Although the *RatioP* algorithm implements the idea of speeding up its operation by selecting a new point by dividing the segment at a given ratio, which ensures its proximity to the point of the current minimum, it otherwise remains practically passive since it does not consider the data accumulated during the operation. In addition, the algorithm is a first-order method, so one cannot expect high speed from it.

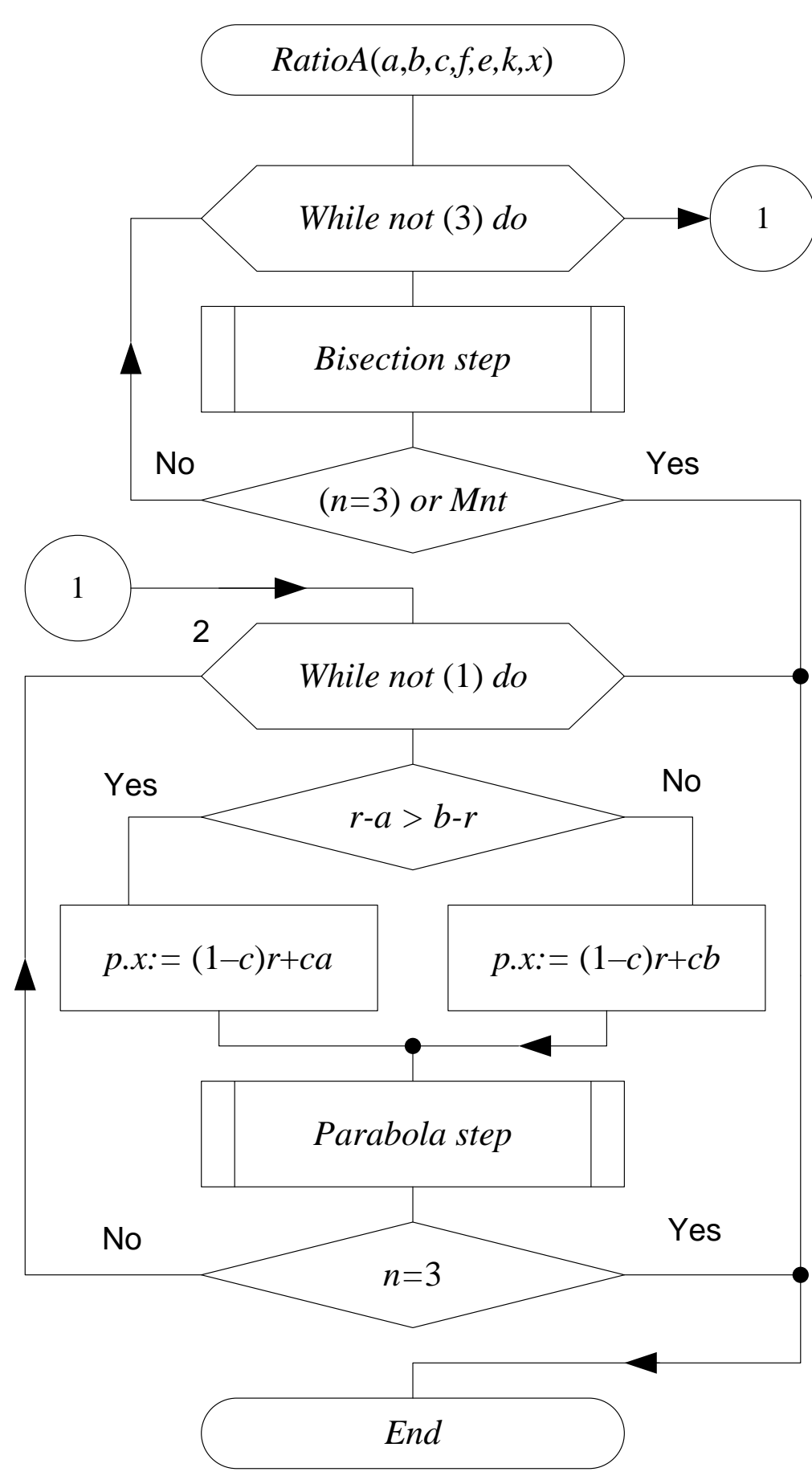

Fig. 6. Block diagram of the *RatioA* algorithm

Since the algorithm quickly recognizes monotone functions and functions with a flat bottom, its efficiency is determined by the speed of working with strictly unimodal functions that have a single minimum. It can detect it via faster methods on the basis of the data obtained during the algorithm's operation. One of them is the parabola method, which is a fast second-order method. It would be possible to apply well-known third-order methods, but they are applicable only to smooth functions since they use their derivative [12].

The flow chart of the *RatioA* algorithm is shown in Fig. 6.

To implement this idea, we first need to find a point with a smaller ordinate, which is in the interval between two previously evaluated points. The problem can then use the parabola method, which will yield a guaranteed result.

To detect the first such point, we can use the *RatioP* method with the parameter $c = 0.5$, working in bisection mode. Using $c < 0.5$ will slow down the algorithm if such a point is close to the edge of the segment.

Having found such a point, we can find the abscissa of the parabola extreme via the following formula [21]:

$$r = \frac{w_{i-1}.y\left(w_i^2.x - w_{i+1}^2.x\right) + w_i.y\left(w_{i+1}^2.x - w_{i-1}^2.x\right) + w_{i+1}.y\left(w_{i-1}^2.x - w_i^2.x\right)}{w_{i-1}.y\left(w_i.x - w_{i+1}.x\right) + w_i.y\left(w_{i+1}.x - w_{i-1}.x\right) + w_{i+1}.y\left(w_{i-1}.x - w_i.x\right)},$$

where $w_{i-1}, w_i, w_{i+1}$ are the three mentioned two-component points of the form $w\ (x, y)$, which satisfy the conditions

$$w_{i-1}.x < w_i.x < w_{i+1}.x, \quad w_{i-1}.y > w_i.y < w_{i+1}.y. \tag{2}$$

The first of them is mandatory. If two adjacent points have the same ordinates, then a flat bottom is detected, and the process must stop. The abscissa of any of these points will be the problem solution.

The parameter $j$ was used to set the ratio $c = 10^{j/2}$. Table 8 shows that the minimum $\Sigma k = 196$ occurs when $c = 0.001$. With this ratio $c$ value, the algorithm has the highest speed.

Table 8. Total number Σk of function evaluations when minimizing them for different values of the parameter *j*

| *j* | –15 | –14 | –13 | –12 | –11 | –10 | –9 | –8 | –7 | –6 | –5 | –4 | –3 | –2 |
|---|---|---|---|---|---|---|---|---|---|---|---|---|---|---|
| Σ*k* | 232 | 226 | 218 | 212 | 205 | 203 | 202 | 201 | 199 | 196 | 209 | 222 | 232 | 278 |

Table 9. Evaluations of all functions 1 – 20

| *N* | *Bisec* | *Gold* | *RatioP* | *RatioA* |
|---|---|---|---|---|
| 1 | 36 | 26 | 3 | 3 |
| 2 | 32 | 25 | 6 | 6 |
| 3 | 36 | 28 | 6 | 6 |
| 4 | 36 | 27 | 4 | 4 |
| 5 | 36 | 29 | 10 | 8 |
| 6 | 28 | 21 | 4 | 4 |
| 7 | 36 | 28 | 23 | 20 |
| 8 | 38 | 29 | 21 | 27 |
| 9 | 40 | 30 | 27 | 21 |
| 10 | 34 | 27 | 18 | 11 |
| 11 | 44 | 33 | 26 | 12 |
| 12 | 36 | 28 | 20 | 4 |
| 13 | 36 | 27 | 18 | 14 |
| 14 | 70 | 51 | 21 | 13 |
| 15 | 40 | 30 | 26 | 13 |
| 16 | 38 | 29 | 22 | 11 |
| 17 | 42 | 32 | 23 | 16 |
| 18 | 38 | 29 | 20 | 10 |
| 19 | 38 | 30 | 22 | 13 |
| 20 | 38 | 28 | 21 | 11 |
| Σ*k* | 772 | 587 | 341 | 227 |
| Relat | 3.40 | 2.59 | 1.50 | 1 |

Table 10. Evaluations of strictly unimodal functions 7 – 20

| *N* | *Bisec* | *Gold* | *RatioP* | *RatioA* |
|---|---|---|---|---|
| 1 | – | – | – | – |
| 2 | – | – | – | – |
| 3 | – | – | – | – |
| 4 | – | – | – | – |
| 5 | – | – | – | – |
| 6 | – | – | – | – |
| 7 | 36 | 28 | 23 | 20 |
| 8 | 38 | 29 | 21 | 27 |
| 9 | 40 | 30 | 27 | 21 |
| 10 | 34 | 27 | 18 | 11 |
| 11 | 44 | 33 | 26 | 12 |
| 12 | 36 | 28 | 20 | 4 |
| 13 | 36 | 27 | 18 | 14 |
| 14 | 70 | 51 | 21 | 13 |
| 15 | 40 | 30 | 26 | 13 |
| 16 | 38 | 29 | 22 | 11 |
| 17 | 42 | 32 | 23 | 16 |
| 18 | 38 | 29 | 20 | 10 |
| 19 | 38 | 30 | 22 | 13 |
| 20 | 38 | 28 | 21 | 11 |
| Σ*k* | 568 | 431 | 308 | 196 |
| Relat | 2.90 | 2.20 | 1.57 | 1 |

Tables 9 and 10 show comparative data for the four algorithms. The fastest algorithm was *RatioA*. The data presented in Table 10 show that functions 7–20 work approximately 1.5 times faster than *RatioP*, 2 times faster than GSS, and at least 2.5 times faster than the bisection search. It is clear that the *RatioA* algorithm works even faster.

## 9. Modernization of Brent's method

The fastest method for minimizing unimodal functions is the combined Brent's method [9]. The method combines the speed of the parabola method and the reliability of the golden section. *The RatioP* algorithm works faster than GSS does; thus, the golden section in Brent's algorithm can be replaced with it in the hope that this will yield a certain positive result. Another idea for improving Brent's method is to expand its capabilities owing to its ability to recognize strictly monotone functions and functions with a flat bottom.

When modernizing Brent's method, a program code written in ALGOL-60 [10] is used as a prototype, in which an errata is taken into account [11]. As a result, a code in Delphi was obtained. Tables 11 and 12 show the comparative results of the modernized *BrentM* code.

As shown in Table 11, the *BrentM* algorithm works faster for all types of functions: 1.69 times faster than Brent's method (*Brent* title), almost 4 times faster than the bisection search, and almost 3 times faster than the GSS. Interestingly, even the performance of the passive algorithm *RatioP* is practically not inferior to that of Brent's algorithm—1.67 versus 1.69. Compared with Brent's algorithm, the active algorithm *RatioA* is more efficient, as it works almost 1.5 times faster. In a certain sense, the advantages of the *BrentM* algorithm are that it works significantly faster than Brent's algorithm when minimizing functions 1–6 in Table 3. Therefore, comparative evaluation data for functions 7–20 in Table 12, when all methods operate under the same conditions, are of interest.

As shown in Table 11, the modernizing *BrentM* algorithm works faster for all types of functions in Table 3, including 1.69 times faster than Brent's algorithm, almost 4 times faster than the bisection search, and almost 3 times faster than the GSS. Interestingly, the *RatioP* algorithm is almost as fast as the Brent's algorithm—1.67 versus 1.69. Compared with Brent's algorithm, the active algorithm *RatioA* is more efficient, as it works almost 1.5 times faster.

In a certain sense, the advantages of the *BrentM* algorithm are that it works significantly faster than Brent's algorithm when minimizing functions 1–6 of Table 3. Therefore, comparative evaluation data for functions 7–20 in this table, when all methods operate under the same conditions, are of interest.

Table 12 shows the results of minimizing these functions. It is clear that the *BrentM* works faster for them as well, although the comparative indicators are not as impressive. Thus, the algorithm works 1.22 times faster than Brent's algorithm, approximately 3 times faster than bisection search, and approximately 2.5 times faster than GSS.

Table 11. Evaluations of all Functions 1 – 20

| *N* | *RatioP* | *RatioA* | *Brent* | *BrentM* |
|---|---|---|---|---|
| 1 | 3 | 3 | 22 | 3 |
| 2 | 6 | 6 | 22 | 6 |
| 3 | 6 | 6 | 25 | 6 |
| 4 | 4 | 4 | 22 | 3 |
| 5 | 10 | 8 | 24 | 6 |
| 6 | 4 | 4 | 16 | 4 |
| 7 | 23 | 20 | 26 | 23 |
| 8 | 21 | 27 | 27 | 22 |
| 9 | 27 | 21 | 21 | 6 |
| 10 | 18 | 11 | 12 | 13 |
| 11 | 26 | 12 | 12 | 12 |
| 12 | 20 | 4 | 6 | 4 |
| 13 | 18 | 14 | 10 | 10 |
| 14 | 21 | 13 | 39 | 12 |
| 15 | 26 | 13 | 13 | 18 |
| 16 | 22 | 11 | 9 | 12 |
| 17 | 23 | 16 | 11 | 11 |
| 18 | 20 | 10 | 9 | 9 |
| 19 | 22 | 13 | 10 | 13 |

Table 12. Evaluations of strictly unimodal functions 7 – 20

| *N* | *RatioP* | *RatioA* | *Brent* | *BrentM* |
|---|---|---|---|---|
| 1 | – | – | – | – |
| 2 | – | – | – | – |
| 3 | – | – | – | – |
| 4 | – | – | – | – |
| 5 | – | – | – | – |
| 6 | – | – | – | – |
| 7 | 23 | 20 | 26 | 23 |
| 8 | 21 | 27 | 27 | 22 |
| 9 | 27 | 21 | 21 | 6 |
| 10 | 18 | 11 | 12 | 13 |
| 11 | 26 | 12 | 12 | 12 |
| 12 | 20 | 4 | 6 | 4 |
| 13 | 18 | 14 | 10 | 10 |
| 14 | 21 | 13 | 39 | 12 |
| 15 | 26 | 13 | 13 | 18 |
| 16 | 22 | 11 | 9 | 12 |
| 17 | 23 | 16 | 11 | 11 |
| 18 | 20 | 10 | 9 | 9 |
| 19 | 22 | 13 | 10 | 13 |

| 20 | 21 | 11 | 9 | 11 |
|---|---|---|---|---|
| Σ$k$ | 341 | 227 | 345 | 204 |
| Relat1 | 1 | 0.67 | 1.01 | 0.60 |
| Relat2 | 1.50 | 1 | 1.52 | 0.90 |
| Relat3 | 0.99 | 0.66 | 1 | 0.59 |
| Relat4 | 1.67 | 1.11 | 1.69 | 1 |

| 20 | 21 | 11 | 9 | 11 |
|---|---|---|---|---|
| Σ$k$ | 308 | 196 | 214 | 176 |
| Relat1 | 1 | 0.64 | 0.69 | 0.57 |
| Relat2 | 1.57 | 1 | 1.09 | 0.90 |
| Relat3 | 1.44 | 0.92 | 1 | 0.82 |
| Relat4 | 1.75 | 1.11 | 1.22 | 1 |

During the computations, Brent's algorithm sometimes produces an error. This occurs when monotone functions are minimized. The error fits into $3\varepsilon_0$, which was warned about by the author of the method [9]. In the modernized Brent’s algorithm, such an error is excluded because the *BrentM* procedure is able to quickly recognize monotone functions. For this class of functions, this algorithm always gives an exact result.

## 10. Conclusion

This paper proposes a new method of ratio section search of the uncertainty interval, which is designed to minimize unimodal functions. The method is able to recognize monotone functions and functions with a flat bottom, which helps to increase its performance, as measured by the number of minimized function evaluations. On the basis of the data obtained by minimizing twenty unimodal functions of various types, a comparison of the performance of the developed method with that of the classical bisection search GSS search was performed. For all types of functions in the bisection mode, the method works faster than bisection by 1.65 times and the GSS by 1.26 times. After optimization, the algorithm with the ratio $c$ begins to work faster by 2.26 and 1.72 times, respectively. After the strictly unimodal functions of Table 3 were minimized, these criteria were 1.84 and 1.4.

The combined Brent’s method, which is considered the fastest method for minimizing unimodal functions, has also been modernized. After the golden section procedure is replaced with the ratio section procedure, a comparison of the performance indicators of the classical and modernized Brent’s methods is made on the basis of the analysis of the computed data. The modernized method works 1.69 times faster than the classical method does. Moreover, the modernized Brent’s method works approximately 4 times faster than the classical bisection search and approximately 3 times faster than the GSS.

**Data Availability**
Data available on request from the authors.

**Ethical Approval**
Author not conducted any human or animal experiments for this article.

**Consent to Participate**
Not needed.

**Consent for Publication**
Not needed.

**Conflict of Interest**
The author declared no competing interests.

**Funding**
No funding.

**Author information**
Siberian Federal University, Polytechnic Institute, Kirensky str. 26, 66074 Krasnoyarsk, Russian Federation.

**References**

[1] Vieira, D.A.G., Lisboa, A. C. ”Line search methods with guaranteed asymptotical convergence to an improving local optimum of multimodal functions”, European Journal of Operational Research, Vol.235, Issue 1, 38 – 46, (2014). https://doi.org/10.1016/j.ejor.2013.12.041.
[2] Waeber, R., Frazier, P. I., Henderson, S. G.”Bisection search with noisy responses”. SIAM Journal on Control and Optimization, 51.3 2261 – 2279, (2013). https://doi.org/10.1137/120861898
[3] Orseau, L., Hutter, M. Line Search for Convex Minimization”. arXiv preprint arXiv:2307.16560. (2023). https://doi.org/10.48550/arXiv.2307.16560

[4] Sharma,P. ”Bisection method or binary-search method role and purposes”. Pranjana: The Journal of Management Awareness, 26 (1 and 2), 107 – 112, (2023). https://doi.org/10.5958/0974-0945.2023.00010.4
[5] Kiefer, J.”Sequential minimax search for a maxim`um”, Proceedings of the American Mathematical Society, 4 (3) 502 – 506. (1953). https://doi.org/10.2307/2032161
[6] Pejic,D., Arsic, M. ”Minimization and Maximization of Functions: Golden-Section Search in One Dimension”. In: Milutinovic, V., Kotlar, M. (eds) Exploring the DataFlow Supercomputing Paradigm. Computer Communications and Networks. Springer, Cham (2019). https://doi.org/10.1007/978-3-030-13803-5_3
[7] Abdelouahhab, M., Manar, S., Benhida, R.”Optimization of the Reaction Temperature During Phosphoric Acid Production Using Fibonacci Numbers and Golden Section Methods”. Chemistry Africa 7, 4017 – 4029, (2024). https://doi.org/10.1007/s42250-024-00971-w
[8] Agushaka, J. O., Ezugwu, A. E., Abualigah, L.”Gazelle optimization algorithm: a novel nature-inspired metaheuristic optimizer”. Neural Comput & Applic 35, 4099 – 4131 (2023). https://doi.org/10.1007/s00521-022-07854-6
[9] Brent, R. P.”Algorithms For Minimization Without Derivatives”. Mathematics of Computation, 19(5), (2002). https://doi.org/10.2307/2005713
[10] Dekker, T. J.”Finding a zero by means of successive linear interpolation”, in Dejon, B.; Henrici, P., Constructive Aspects of the Fundamental Theorem of Algebra, London: Wiley-Interscience, (2009). ISBN 978-0-471-20300-1
[11] Gegenfurtner,K. R. ”PRAXIS: Brent's algorithm for function minimization”. Behavior Research Methods, Instruments, & Computers 24, 560 – 564 (1992). https://doi.org/10.3758/BF03203605
[12] Kodnyanko, V. A., Grigorieva, O. A., Strok, L.V.”Combined Newton's third-order convergence method for minimize one variable functions”. Radio Electronics, Computer Science, Control, (2), 48 – 55 (2021). https://doi.org/10.15588/1607-3274-2021-2-5
[13] Steffen, V., Della Pasqua, C. C., De Oliveira, M. S., Da Silva, E. A. ”Halving interval guaranteed for Dekker and Brent root finding methods”. Examples and Counterexamples, 100173. (2024) https://doi.org/10.2139/ssrn.4702916
[14] Schommer, A., Xavier, M. A., Morrey, D., Collier, G.”Real-Time Deployment Strategies for State of Power Estimation Algorithms” No. 2024-01-2198). SAE Technical Paper (2024). https://doi.org/10.4271/2024-01-2198
[15] Rao, S. S.”Engineering optimization: theory and practice”. John Wiley & Sons, Hoboken, New Jersey. (2009).
[16] Wilde, D. J. Optimum seeking methods. Englewood Cliffs, NJ: Prentice Hall, Inc. 423 p. (1964),
[17] V. Beiranvand, W. Hare,Y. Lucet. ”Best practices for comparing optimization algorithms”. Optim. Eng. 18, 815 – 848 (2017).
https://doi.org/10.1007/s11081-017-9366-1
[18] Kodnyanko, V. A.”Economical dichotomous search for minimizing one-variable functions”. Radio Electronics, Computer Science, Control, (3), 34 – 39, (2019). https://doi.org/10.15588/1607-3274-2019-3-4
[19] Teti, D. Delphi Cookbook. Second Edition. Packt Publishing Ltd., 2016. ISBN 978-1-78528-742-8.
[20] Rao, C. R.”Toutenburg Statistics (3rd ed.) ”. Berlin: Springer. (2008). ISBN 978-3-540-74226-5. 1998.